\documentclass[11pt]{amsart}
\usepackage{amsmath,amsthm,amssymb,amstext,amsxtra}
\usepackage[all]{xy}
\usepackage[alphabetic,backrefs,lite]{amsrefs}
\usepackage{mathrsfs}  
\usepackage{mathtools} 
\usepackage{colonequals}
\usepackage{enumitem}

\usepackage{graphicx}
\usepackage{comment}

\usepackage[backref, colorlinks=true, linkcolor=blue, citecolor=blue]{hyperref}
\usepackage[]{xcolor}
\usepackage[top=1in, bottom =1in, left=1in, right = 1in]{geometry}

\newcommand{\Z}{\mathbb{Z}}
\newcommand{\Q}{\mathbb{Q}}

\newcommand{\A}{\mathbb{A}}

\newcommand{\F}{\mathbb{F}}

\renewcommand{\P}{\mathbb{P}}

\DeclareMathOperator{\Br}{Br}

\newcommand{\Div}{\mathrm{Div}}

\DeclareMathOperator \disc {disc}
\DeclareMathOperator{\bPic}{{\bf Pic}}
\DeclareMathOperator{\Gr}{Gr}

\usepackage[OT2,T1]{fontenc}
\DeclareSymbolFont{cyrletters}{OT2}{wncyr}{m}{n}
\DeclareMathSymbol{\Sha}{\mathalpha}{cyrletters}{"58}

\theoremstyle{plain}
\newtheorem{theorem}[equation]{Theorem}

\newtheorem{corollary}[equation]{Corollary}

\theoremstyle{definition}

\theoremstyle{remark}
\newtheorem{remark}[equation]{Remark}

\title{A threefold violating a local-to-global principle for rationality}
\author{Sarah Frei}
\address{Department of Mathematics, Dartmouth College, 27 N. Main Street, Hanover, NH 03755}
\email{sarah.frei@dartmouth.edu}
\urladdr{https://math.dartmouth.edu/\~{}sfrei}

\author{Lena Ji}
\address{Department of Mathematics, University of Michigan, 530 Church Street, Ann Arbor, MI 48109-1043}
\email{lenaji.math@gmail.com}
\urladdr{http://www-personal.umich.edu/\~{}lenaji}

\thanks{L.J. is supported by NSF MSPRF grant DMS-2202444.}

\subjclass[2020]{Primary: 14E08. Secondary: 14G12, 14J20, 11G30.}

\begin{document}

\begin{abstract}
    In this note we construct an example of a smooth projective threefold that is irrational over \(\mathbb Q\) but is rational at all places. Our example is a complete intersection of two quadrics in \(\mathbb P^5\), and we show it has the desired rationality behavior by constructing an explicit element of order \(4\) in the Tate--Shafarevich group of the Jacobian of an associated genus \(2\) curve.
\end{abstract}

\maketitle

\section{Introduction}
Let \(X\) be a smooth projective variety over a global field \(k\), and suppose that for every place \(v\) of \(k\), the variety \(X\times_k k_v\) is \(k_v\)-rational (i.e. birational to projective space over \(k_v\)). Is \(X\) \(k\)-rational?

In dimension \(1\), \(k_v\)-rationality at all places implies \(k\)-rationality, since conics satisfy the Hasse principle, and existence of a point characterizes rationality for conics. In dimension \(2\), this local-to-global principle for rationality again holds if \(X\) is a del Pezzo surface of degree at least \(5\) \cite{Manin}*{Section IV}. In smaller degree, however, there are examples of del Pezzo surfaces that are \(\mathbb Q\)-unirational and rational over all completions of \(\mathbb Q\), but have a Brauer group obstruction to \(\mathbb Q\)-rationality \cite{Coray-Tsfasman}*{Example 3.3 and the following remark}.

We study the case of threefolds, constructing a threefold complete intersection of two quadrics \(X\subset \mathbb P^5\) over \(\mathbb Q\) that is rational at all places but \emph{irrational} over \(\mathbb Q\). Note that \(\Br X=\Br\mathbb Q\) for such \(X\) \cite{CTSSD}, so the rationality obstruction used in the surface example of \cite{Coray-Tsfasman} vanishes.
We are not aware of any previous such explicit examples of dimension \(\geq 3\) in the literature.

\begin{theorem}\label{thm:main-example} Let \(\mathcal X=Q_1\cap Q_2\subset \mathbb P^5_{\mathbb Z, [u:v:w:x:y:z]}\), where the quadrics \(Q_1,Q_2\) are defined by
\begin{equation*}
    \begin{split}
        Q_1 &= uv + uw - 4vw + 2vz + 2wz + x^2 - 2xz + y^2 - z^2, \\
        Q_2 &= uv - uw + uy - 2v^2 + 2vx - 2wy + 2wz + 2xz.
    \end{split}
\end{equation*}

Then \(X\coloneqq\mathcal X\times_{\mathbb Z}\mathbb Q\) is a smooth projective threefold that is \(\mathbb Q\)-unirational, \(\mathbb R\)-rational, and \(\mathbb Q_p\)-rational for all primes \(p\), but is \emph{irrational} over \(\mathbb Q\). Furthermore, the reduction of \(\mathcal X\) modulo \(p\) is \(\mathbb F_p\)-rational for all primes \(p\neq 2\).
\end{theorem}
The intermediate Jacobian of \(X\) is the Jacobian of the curve \(C : z^2=-t^6 - 3t^5 + 2t^4 + 3t^3 - 3t^2 - 3t - 2\), which appeared in \cite{Bruin-Stoll-experiment}. 
Bruin--Stoll, when combined with recent work of Fisher--Yan \cite{FisherYan}, show that \(\bPic^1_C\) is a nontrivial element of \(\Sha(\bPic^0_C)\).
The recent rationality criterion of Hassett--Tschinkel over \(\mathbb R\) \cite{HT-intersection-quadrics} and Benoist--Wittenberg over arbitrary fields \cite{BW-IJ} (see also~\cite{HT-cycle-class-maps} for \(k\subset \mathbb C\)) shows that the existence of a \(k\)-point on a certain 2-cover of \(\bPic^1_C\) determines rationality of \(X\) (see Section~\ref{sec:background}); thus, we show Theorem~\ref{thm:main-example} by constructing an explicit example where this \(2\)-cover violates the Hasse principle (and, more precisely, has order \(4\) in the Tate--Shafarevich group).

In \cite{FJSVV} the authors, together with Sankar, Viray, and Vogt, construct examples of conic bundle threefolds over \(\mathbb Q\) that are irrational over \(\mathbb R\) (and hence irrational over \(\mathbb Q\)), rational over \(\mathbb C\), and become rational modulo all primes of good reduction (for the discriminant double cover). Earlier, examples of threefold intersections of two quadrics that are irrational over \(\mathbb R\) and rational modulo all primes of good reduction appeared implicitly in Hassett--Tschinkel's work \cite{HT-intersection-quadrics}*{Construction 1, Theorem 36}. The difficulty in constructing an intersection of quadrics \(X\) as in Theorem~\ref{thm:main-example} lies in distinguishing between rationality over \(\mathbb Q\) and \(\mathbb R\), and in determining the behavior at the bad primes. Note that \(X\) cannot have everywhere good reduction, since otherwise its intermediate Jacobian would be a nontrivial abelian variety with everywhere good reduction, which is impossible by \cites{Abraskin,Fontaine}.

The Hasse principle for smooth intersections of two quadrics in \(\mathbb P^5\) is an open question. Wittenberg showed it holds under the assumptions of Schinzel's hypothesis and finiteness of \(\Sha\) for elliptic curves \cite{Wittenberg-intersection-quadrics}. Recently, Iyer--Parimala showed that if \(X\) contains a \(k_v\)-line for all places \(v\) of \(k\) and if \(\disc(Q_1)=1\), then \(X\) contains a \(k\)-point \cite{IP-period-index}*{Theorem 0.2}. (By \cite{IP-period-index}*{Corollary 10.2}, see also \cite[Corollaire~8.7]{CT-intersection-quadrics}, the condition \(\disc(Q_1)=1\) can be replaced with assuming that \(\Div^1(C)(k)\neq\emptyset\)\footnote{\(\Div^i_C(k)\) denotes the set of \(k\)-rational divisors of degree \(i\). \(\bPic^i_C\) denotes the degree \(i\) component of the Picard scheme, and its \(k\)-points are the \(k\)-rational divisor classes of degree \(i\).}, where \(C\) is the genus \(2\) curve associated to the discriminant of the pencil of quadrics spanned by \(Q_1\) and \(Q_2\).)

A threefold intersection of quadrics defines a Brauer class on the associated genus \(2\) curve \(C\) by taking the Azuyama algebra of the even Clifford algebra associated to the pencil of quadrics (see Section~\ref{sec:background}). If \(k\) is a number field and \(\Div^1(C)(k)\neq\emptyset\), Iyer--Parimala proved the period-index conjecture for the \(2\)-torsion of \(\Sha(\Br C)\coloneqq\ker\left( H^2(C,\mathbb G_m) \to \prod_v H^2(C_v, \mathbb G_m) \right)\) \cite{IP-period-index}*{Corollary 5.4}. As a consequence of Theorem~\ref{thm:main-example} we obtain the following example of a nontrivial Brauer class on a genus \(2\) curve without a degree \(1\) zero-cycle:

\begin{corollary}\label{cor:Sha-Br}
    Let \(Q_1,Q_2\), and \(C\) be as in Theorem~\ref{thm:main-example}. Then \(\Div^1(C)(\Q)=\emptyset\), and the pencil of quadric fourfolds spanned by \(Q_1\) and \(Q_2\) defines a nontrivial \(2\)-torsion class \(\beta\in\Sha(\Br C)\).
\end{corollary}

After the first version of this article appeared on the arXiv, Kunyavskii \cite{Kunyavskii} has used different methods to construct new examples of varieties violating the local-to-global principle for rationality and where the Brauer obstruction vanishes. He produces examples of algebraic torii over any global field, and surface examples over any global field of characteristic \(\neq 2\).

\subsection*{Acknowledgments} We are grateful to Anthony V\'arilly-Alvarado for computing advice, to Tom Fisher for suggesting a simplifying change of coordinates, and to Tom Fisher and Jiali Yan for correspondence about \cite{FY-code, FisherYan}. We also thank Eran Assaf, Asher Auel, Nils Bruin, Sebastian Casalaina-Martin, Brendan Creutz, Sam Frengley, Jack Petok, and Bianca Viray for helpful conversations. The computations in Section~\ref{sec:construction} were done in \texttt{Magma} \cite{Magma}.

\section{Intersections of two quadrics in \(\mathbb P^{5}\) and genus \(2\) curves}\label{sec:background}

In this section we recall results about the geometry of a pencil of two quadrics. The geometry of the variety of maximal linear spaces in an intersection of two quadrics is a rich theory that has been widely studied. We will only address the case of quadrics in \(\mathbb P^5\) over fields of characteristic not \(2\), because this is the generality that we will require; we refer the reader to \cite{Reid-thesis, CTSSD} for other dimensions and to \cite{BW-IJ}*{Section 4} for arbitrary characteristic.

Let \(k\) be a field of characteristic not equal to \(2\), and let \(X=Q_1\cap Q_2\subset \mathbb P^5_k\) be a complete intersection of two quadrics. Let \(M_1\) and \(M_2\) denote the Gram matrices of \(Q_1\) and \(Q_2\) with respect to the same basis, and define the polynomial \[f(t)= -\det(M_1 - t M_2).\] Then \(X\) is a smooth threefold if and only if the polynomial \(f(t)\) is not identically zero and has \(6\) distinct roots \cite{Reid-thesis}*{Proposition 2.1}.

If \(X\) is smooth, then its intermediate Jacobian is the Jacobian of the genus \(2\) curve \(C\) defined by \(z^2=f(t)\), and its Fano variety of lines \(F_1(X)\) is a smooth surface \cite{Reid-thesis}*{Theorem 2.6}. It is classically known that \(X\) is \(k\)-rational if \(F_1(X)(k)\neq\emptyset\), since projection from a line on \(X\) gives a rationality construction \cite{CTSSD}*{Proposition 2.2}. Wang showed that \(F_1(X)\) is a torsor under \(\bPic^0_C\), and that \(2[F_1(X)] = [\bPic^1_C]\) as torsors under \(\bPic^0_C\) over \(k\) \cite{Wang-maximal-linear-spaces}*{Theorem 1.1}. By extending the Clemens--Griffiths intermediate Jacobian rationality obstruction to a torsor condition over non-closed fields and using Wang's work, Hassett--Tschinkel (over \(\mathbb R\)) and Benoist--Wittenberg (over arbitrary fields) characterized rationality of \(X\) over any field:
\begin{theorem}[\citelist{\cite{HT-intersection-quadrics}*{Theorem 36} \cite{BW-IJ}*{Theorem A} \cite{HT-cycle-class-maps}*{Theorem 24}}]\label{thm:rationality-intersection-quadrics}
    Let \(X\subset\mathbb P^5_k\) be a smooth complete intersection of two quadrics. Then \(X\) is \(k\)-rational if and only if \(F_1(X)(k)\neq\emptyset\).
\end{theorem}
In particular, \(X\) is rational at all places but irrational over \(\mathbb Q\) if and only if the surface \(F_1(X)\) violates the Hasse principle; this failure of the Hasse principle should be explained by a Brauer--Manin obstruction on \(F_1(X)\) \cite{Skorobogatov-torsors}*{Theorem 6.2.3}.

The complete intersection \(X\) defines a natural Brauer class on the genus \(2\) curve \(C\) as follows. Let \(\mathcal Q\to\mathbb P^1\) be the pencil of quadric fourfolds spanned by \(Q_1\) and \(Q_2\). The even Clifford algebra of this pencil defines an Azumaya algebra on the discriminant double cover \(C\to\mathbb P^1\) and hence a Brauer class \(\beta\in\Br C\). The proof of \cite[Lemma 5.2.3]{ABB-quadrics} shows:
\begin{theorem}[{c.f.~\cite[Lemma 5.2.3]{ABB-quadrics}}]\label{thm:line-BrC}
Let \(k\) be a field of characteristic \(\neq 2\), and let \(X\subset\mathbb P^5_k\) be a smooth complete intersection of two quadrics with \(X(k)\neq\emptyset\). Then \(F_1(X)(k)\neq\emptyset\) if and only if \(\beta\in \Br C\) is trivial.
\end{theorem}

As explained above, a pencil of quadrics naturally gives rise to a genus \(2\) curve. In the other direction, given a genus \(2\) curve \(C\) over \(k\), Bhargava--Gross--Wang characterize when \(C\) can arise from an intersection of two quadrics in the above manner \cite{BGW}*{Theorem 24}. Furthermore, in this case, they give conditions for the existence of lines on \(X\).
\begin{theorem}[\cite{BGW}*{Theorem 29}]\label{thm:soluble-orbits}
Let \(C: z^2=f(t)\) be a smooth projective curve of genus \(2\) defined over a field \(k\) of characteristic \(\neq 2\). The following two conditions are equivalent:
\begin{enumerate}
    \item There exist quadrics \(Q_i\) over \(k\) such that \(f(t)=-\det(M_1 - t M_2)\) and \(F_1(Q_1\cap Q_2)(k)\neq\emptyset\).
    \item \(\Div^1(C)(k)\neq\emptyset\).
\end{enumerate}
\end{theorem}

The above result does \emph{not} imply that every \(C\) over a global field \(k\) with local points everywhere is obtained from an intersection of quadrics \(X\subset\mathbb P^5\) such that \(X\) has \(k_v\)-lines everywhere: it may be the case that no single choice of quadrics over \(k\) simultaneously works for all places.
 
 Indeed, the existence of such an \(X\) implies the \(\bPic^0_C\)-torsor \(F_1(X)\) is locally soluble and satisfies \(2[F_1(X)]=[\bPic^1_C]\); hence, \(F_1(X)\) is a locally soluble \(2\)-cover of \(\bPic^1_C\) (see \cite{BGW}*{page 2} for definitions), and in particular the class of \(\bPic^1_C\) is divisible by \(2\) in the Tate--Shafarevich group of \(\bPic^0_C\).
 Bhargava--Gross--Wang showed that this necessary condition is sufficient:
\begin{theorem}[\cite{BGW}*{Theorem 31}]\label{thm:locally-soluble-orbits}
Let \(k\) be a global field of characteristic \(\neq 2\), and let \(C: z^2=f(t)\) be a smooth projective curve of genus \(2\) defined over \(k\) such that \(\Div^1(C)(k_v)\neq\emptyset\) for all places \(v\) of \(k\). The following two conditions are equivalent:
\begin{enumerate}
    \item There exist quadrics \(Q_i\) defined over \(k\) such that \(f(t)=-\det(M_1 - t M_2)\) and \(F_1(Q_1\cap Q_2)(k_v)\neq\emptyset\) for all places \(v\).
    \item \(\bPic^1_C\) admits a locally soluble \(2\)-cover over \(k\).
\end{enumerate}
\end{theorem}

Combining Theorems~\ref{thm:rationality-intersection-quadrics}, \ref{thm:soluble-orbits}, and \ref{thm:locally-soluble-orbits}, together with the fact that \(2[F_1(X)]=[\bPic^1_C]\) as \(\bPic^0_C\)-torsors, gives the following sufficient condition for the existence of an irrational \(X\) that is everywhere locally rational.
\begin{corollary}\label{cor:criterion}
Let \(k\) be a global field of characteristic \(\neq 2\). Let \(C: z^2=f(t)\) be a smooth projective curve of genus 2 defined over k. Assume the following conditions hold:
\begin{enumerate}
    \item\label{item:Pic^1-no-Q-point} \(\bPic^1_C(k) = \emptyset\), and
    \item\label{item:4-torsion-Sel} \(\bPic^1_C\) admits a locally soluble \(2\)-cover over \(k\).
\end{enumerate}
Then there exists a smooth complete intersection of two quadrics \(X \subset \mathbb P^5_k\) whose intermediate Jacobian is \(\bPic^0_C\) and is such that \(X\) is \(k_v\)-rational for all places \(v\) of \(k\), but \(X\) is not \(k\)-rational.
\end{corollary}
In particular, conditions~\eqref{item:Pic^1-no-Q-point} and~\eqref{item:4-torsion-Sel} imply that \(\Div^1(C)\) violates the Hasse principle, and that \(\Sha(\bPic^0_C)\) contains an element of order \(4\).

Note that even if \(C\) satisfies the conditions in Corollary~\ref{cor:criterion}, it need not be the case that {\it every} complete intersection of quadrics with associated genus \(2\) curve \(C\) is rational at all places. (In the language of \cite{BGW}, the existence of locally soluble orbits does not imply that every orbit is locally soluble; see \cite{BGW}*{Proof of Theorem 31}.)

\begin{remark}
    Condition~\eqref{item:Pic^1-no-Q-point} in Corollary~\ref{cor:criterion} is not necessary for the
    existence of a threefold complete intersection of two quadrics violating the local-to-global principle for rationality. Indeed, if \(C\) has a \(k\)-point and if the \(2\)-Selmer group of \(\bPic^0_{C}\) is nontrivial, then such threefolds with associated \(C\) exist \cite[Theorem 11]{ShankarWang}. However, if \(\bPic^1_C(k)\neq\emptyset\), we do not know a way in practice to determine, from explicit equations for \(X\), whether the torsor \(F_1(X)\) is nontrivial.
\end{remark}

\subsection{Genus \(2\) curves with \(\bPic^1_C\) violating the Hasse principle}
We construct our example by searching the literature for genus \(2\) curves satisfying the conditions in Corollary~\ref{cor:criterion}, performing the procedure described in \cite{BGW}*{Section 2} (for suitable choices of \(s\) and \(\alpha\)) to recover a complete intersection of quadrics \(X\) whose discriminant is associated to this curve, and then verifying that \(F_1(X)\) has local points everywhere.

Any genus \(2\) curve \(C\) over a local field \(k_v\) necessarily has \(\bPic^1_C(k_v)\neq\emptyset\) \cite{PoonenStoll}*{Corollary 4, Footnote 10}, but the existence of a \(k_v\)-rational divisor \emph{class} of degree \(1\) does not necessarily imply that \(\Div^1(C)(k_v)\neq\emptyset\).
Moreover, many genus \(2\) curves \(C\) that have been previously shown to have interesting arithmetic do not satisfy the criteria in \cite{BGW}*{Theorem 24} for \(\bPic^0_C\) to arise globally as the intermediate Jacobian of a complete intersection of two quadrics. For example, a smooth genus \(2\) curve \(C\) over \(\mathbb Q\) with no \(\mathbb R\)-points and with \(\Div^1(C)(\mathbb Q_p)\neq\emptyset\) for all primes \(p\) has \(\bPic^1_C(\mathbb Q)=\emptyset\) by \cite{PoonenStoll}*{Theorem 11} (see \cite{AAFH}*{Section 3} for a description of the Brauer--Manin obstruction, and \cite{PoonenStoll}*{\(t>0\) case of Proposition 26} for an example of such a curve).
However, a genus \(2\) curve with \(C(\mathbb R)=\emptyset\) will never arise from an intersection of quadrics over \(\mathbb R\) \cite{BGW}*{Section 7.2}.

The example constructed in \cite{PoonenStoll}*{Proposition 28} with both \(C\) and \(\bPic^1_C\) violating the Hasse principle does arise over \(\mathbb Q\) from an intersection of quadrics and satisfies condition~\eqref{item:Pic^1-no-Q-point} in Corollary~\ref{cor:criterion}; however, Poonen--Stoll show conditionally on the Birch--Swinnerton-Dyer conjecture that \(\Sha(\bPic^0_C)\cong \mathbb Z/2\times\mathbb Z/2\), so condition~\eqref{item:4-torsion-Sel} does not hold. See also \cite{flynn-art22} for additional examples of genus \(2\) curves where \(C\) and \(\bPic^1_C\) violate the Hasse principle.

To prove Theorem~\ref{thm:main-example}, we use a genus \(2\) curve \(C\) from \cite{Bruin-Stoll-experiment}. Bruin--Stoll and Fisher--Yan \cite{FisherYan} exhibit \(38\) examples of curves \(C\) satisfying the conditions of Corollary~\ref{cor:criterion}. More precisely, Bruin--Stoll construct \(42\) curves \(C\) such that:
\begin{enumerate}[label=(\alph*)]
    \item \(C\) is locally soluble,
    \item \(\bPic^1_C\) has a locally soluble \(2\)-cover, and
    \item\label{item:BS-conditional} (conditionally) \(\bPic^1_C(\mathbb Q)=\emptyset\).
\end{enumerate}
  For~\ref{item:BS-conditional}, Bruin--Stoll reduce showing \(\bPic^1_C(\mathbb Q)=\emptyset\) to computing the rank of the Jacobian of \(C\). Then, they assume either the Generalized Riemann Hypothesis or the Birch--Swinnerton-Dyer conjecture for \(C\) in order to conditionally compute the rank for each of these \(42\) curves \cite{Bruin-Stoll-experiment}*{Sections~2 and 3.4}. Recently, Fisher--Yan \cite{FisherYan} developed a new method for computing the Cassels--Tate pairing on the \(2\)-Selmer group of a genus \(2\) Jacobian, which allows them show that the Jacobian has expected rank for all but \(4\) of these curves. (The \(4\) unresolved cases are \#2, \#21, \#23, and \#33 \cite{FY-code}.) In particular, they show that \ref{item:BS-conditional} holds unconditionally for these \(38\) genus \(2\) curves \cite[Section 4.3]{FisherYan}.

\section{Construction of the example}\label{sec:construction}

Let \(Q_1,Q_2\) be the quadratic forms in Theorem~\ref{thm:main-example}, and let \(M_i\) be the symmetric matrix associated to \(Q_i\). Then \[ -\det(M_1 - t M_2) = -t^6 - 3t^5 + 2t^4 + 3t^3 - 3t^2 - 3t - 2.\]
The genus \(2\) curve \(z^2=-\det(M_1 - t M_2)\) is curve \#22 in the \texttt{BSD-data.txt} file \cite{BS-data} accompanying \cite{Bruin-Stoll-experiment}. This and all of the computational claims in the proof below can be verified with \texttt{Magma} code provided in our \texttt{Github} repository \cite{FJ-code}.

\subsection{Equations for \(F_1(X)\)}\label{sec:variety-of-lines}
To show that \(X\) has \(\mathbb Q_p\)-lines for all primes \(p\), we show that the Fano variety of lines \(F_1(X)\) has smooth \(\mathbb F_p\)-points for all primes \(p\) and then apply Hensel's lemma. Such points are automatic for primes of good reduction by Lang's theorem, since in the smooth case \(F_1(X)\) is a \(\bPic^0_C\)-torsor. For the primes of bad reduction, we must work with the explicit equations for \(F_1(X)\). To do this, we work in a standard affine patch of \(\Gr(2,6)\); see for example \cite{Harris}*{Example 6.6}. In the accompanying \texttt{Magma} code, we use the affine patch \(\mathbb A^8_{(t_j)}\) of \(\Gr(2,6)\) parametrizing lines given by the row space of the matrix 
 \[\begin{pmatrix} t_1 & 1 & 0 & t_3 & t_5 & t_7 \\ t_2 & 0 & 1 & t_4 & t_6 & t_8 \end{pmatrix}.\]
Explicitly, on this chart of \(\Gr(2,6)\) we consider lines \(\P^1\hookrightarrow \P^5\) of the form \([r:s] \mapsto [t_1r+t_2s : r: s : t_3r+t_4s : t_5r+t_6s : t_7r+t_8s].\)\footnote{We work with this particular affine patch because, on many other patches, \(F_1(X)\) has no smooth \(\F_2\)-points.}
Requiring a line to be contained in the quadric \(Q_i\) defines conditions on the \(t_j\) (eqnsi in the \texttt{Magma} code); together these equations for \(Q_1\) and for \(Q_2\) define the Fano surface of lines on this affine patch. In the \texttt{Magma} code, this chart is called Fanopatch.

\subsection{Proof of Theorem~\ref{thm:main-example} and Corollary~\ref{cor:Sha-Br}}

First, \([1 : 0 : 0 : 0 : 0 : 0 ]\in \P^5_{\Q}\) exhibits a \(\mathbb Q\)-point on \(X\), so \(X\) is unirational over \(\mathbb Q\) \cite{CTSSD}*{Proposition 2.3}.

Next, we show that \(X\) becomes rational at all places. To do this, we will show that \(F_1(X)\) has a \(\Q_v\)-point for all places of \(\Q\) (including the real place), so that \(X\) has a \(\Q_v\)-line and is thus rational over \(\mathbb Q_v\) by \cite{CTSSD}*{Proposition 2.2}. First, \(C\) has two real Weierstrass points, branched over points in \(\P^1\) whose coordinates are approximately \([ -3.26599: 1]\) and \([-1.13643: 1 ]\). In this case \(F_1(X)\) necessarily has a real point \cite{BGW}*{Section 7.2}. 

Now fix a prime \(p\). If \(p\) does not divide the discriminant of \(C\), then the reduction of \(X\) mod \(p\) is smooth \cite[Proposition 2.1]{Reid-thesis}, so \(F_1(X)\) has a \(\mathbb Q_p\)-point (see Section~\ref{sec:variety-of-lines} above). 

Then it remains to check for the primes \(p\) dividing the discriminant of \(C\) that the reduction of \(F_1(X)\) modulo \(p\) contains a smooth \(\F_p\)-point, which we do by checking directly on the affine open chart described in Section~\ref{sec:variety-of-lines}. The two primes dividing the discriminant are \(2\) and \(149743897\).
\begin{itemize}
\item For \(p=2\), a smooth \(\mathbb F_2\)-point is given by \((1, 1, 0, 0, 1, 1, 0, 0)\in \A^8_{\F_2}\) on the affine patch discussed in Section~\ref{sec:variety-of-lines}.
\item
For \(p=149743897\), on the same affine patch, a smooth \(\mathbb F_{149743897}\)-point is given by 
\[(10276, 859210, 113976451, 113430900, 122036333, 94785567, 35411179,
25838500)\in \A^8_{\F_{149743897}}.\]
\end{itemize}
In each case, Hensel's lemma yields a \(\mathbb Q_p\)-point on \(F_1(X)\) and hence an \(\mathbb Q_p\)-line on \(X\). Thus, we have checked that \(X\) is rational at all places of \(\mathbb Q\).

We next consider \(\mathcal X\) over \(\Z\) as given in Theorem~\ref{thm:main-example}, and show that the mod \(p\) reduction \(X_p:=\mathcal{X}\times_\Z \F_p\) is \(\mathbb F_p\)-rational for all primes \(p\neq 2\). For \(p\not\in\{ 2, 149743897\}\) the reduction \(X_p\) is smooth, so this follows from \cite{HT-intersection-quadrics}*{paragraph after Construction 1}. For \(p=149743897\), since we have shown that \(X_{149743897}\) contains a line, \cite{CTSSD}*{Proposition 2.2} implies that \(X_{149743897}\) is \(\mathbb F_{149743897}\)-rational if it is not a cone. Checking that \(X_{149743897}\) is non-conical is equivalent, by \cite{CTSSD}*{Lemma 1.12}, to checking that the Jacobian matrix of \((Q_1,Q_2)\) does not vanish identically at any point. Indeed, the singular locus consists of a single point \([10925789 : 85737939 : 85378598 : 93099029 : 51694582 : 1]\in \P^5_{\F_{149743897}}\), and the Jacobian matrix has rank \(1\) at this point. Thus, we have verified \(\mathbb F_p\)-rationality of \(X_p\) for all \(p\neq 2\). (The mod \(2\) reduction \(X_2\) is reducible and non-reduced.)

It remains to show the irrationality of \(X\) over \(\mathbb Q\). To show this, we claim that \(F_1(X)\) has order \(4\)  in \(\Sha(\bPic^0_C)\); by Theorem~\ref{thm:rationality-intersection-quadrics} this implies \(X\) is irrational over \(\mathbb Q\). Since \(2[F_1(X)]=[\bPic^1_C]\) in the Weil--Ch\^{a}telet group of \(\bPic^0_C\) \cite{Wang-maximal-linear-spaces}, it suffices to show that \(\bPic^1_{C}\) has no \(\Q\)-points. For this, Bruin--Stoll show that if the Jacobian of \(C\) has expected rank (which is \(0\) in this case), then \(\bPic^1_{C}(\mathbb Q)=\emptyset\) \cite{Bruin-Stoll-experiment}*{Section 2}. They then conditionally compute the rank, assuming the Birch--Swinnerton-Dyer conjecture for \(C\). Recent work of Fisher--Yan \emph{unconditionally} computes the rank \cite[Section 4.3]{FisherYan}; details are provided in the Magma code \texttt{g2ctp\_examples.m} in \cite{FY-code}, where this curve is called Bruin Stoll \#22. Thus, \(\bPic^1_{C}\) has no \(\Q\)-points, and the claim follows.

For Corollary~\ref{cor:Sha-Br}, first note that \(\Div^1(C)(\Q)=\emptyset\) because \(\bPic^1_{C}(\mathbb Q)=\emptyset\). Thus, Corollary~\ref{cor:Sha-Br} follows from Theorem~\ref{thm:main-example} (using Theorems~\ref{thm:rationality-intersection-quadrics} and \ref{thm:line-BrC}).
\qed

\begin{remark}
In Theorem~\ref{thm:main-example}, local solubility of \(F_1(X)\) at \(\mathbb R\) is automatic because the curve \(C\) has two real Weierstrass points \cite{BGW}*{Section 7.2}. In general, if \(C\) has any number of real Weierstrass points, Krasnov has characterized when \(X\) contains a real line, by giving an isotopy classification using the signatures of the matrices \(t_0 M_1- t_1 M_2\) as \((t_0,t_1)\) varies over the unit circle in \(\mathbb R^2\) \cite{krasnov-biquadrics}.
\end{remark}

\bibliographystyle{alpha}
\bibliography{references.bib}

\end{document}